\input amstex
\documentstyle{amsppt}
\input algvar.def
\loadbold
\finaltrue
\onefiletrue

\topmatter
\title \SPINC\ Structures and Scalar Curvature Estimates\endtitle
\author S. Goette, U. Semmelmann\endauthor
\rightheadtext{Scalar Curvature Estimates}
\address  Universit\'e de Paris-Sud, D\'epartement de Math\'ematique,
Laboratoire de Topologie et Dynamique, URA D1169 du CNRS, B\^atiment 425,
91405 Orsay Cedex, France\endaddress
\email Sebastian.Goette\@math.u-psud.fr\endemail
\address Mathematisches Institut, Universit\"at M\"unchen,
Theresienstr.~39, D-80333 M\"unchen, Germany\endaddress
\email semmelma\@rz.mathematik.uni-muenchen.de\endemail
\thanks Both authors were supported by a research fellowship of the~DFG
\endthanks
\abstract
In this note,
we look at estimates for the scalar curvature~$\kappa$
of a Riemannian manifold~$M$
which are related to \spinc\ Dirac operators:
We show that one may not enlarge a K\"ahler metric with positive
Ricci curvature without making~$\kappa$ smaller somewhere on~$M$.
We also give explicit upper bounds for~$\min\kappa$
for arbitrary Riemannian metrics on certain
submanifolds of complex projective space.
In certain cases, these estimates are sharp:
we give examples where equality is obtained.
\endabstract
\endtopmatter

\document
\Introduction
\head 0. Introduction\endhead
There is a relation between the positivity of the scalar curvature~$\kappa$
of a compact Riemannian manifold~$M$ and its topology:
Let~$M$ be an $m$-dimensional, compact, orientable, spin manifold.
If the generalized $\Adach$-genus~$\alpha(M)\in KO^{-m}(\roman{pt})$
does not vanish,
then~$M$ carries no metric of positive scalar curvature
by theorems of Lichnerowicz and Hitchin.
The proof uses Dirac operators,
combining the Bochner-Weitzenb\"ock-Lichnerowicz formula (short: BLW-formula)
and the Atiyah-Singer index theorem.
The converse is much harder to establish:
If~$m\ge 5$ and~$M$ is compact,
simply-connected and either not spin or spin with~$\alpha(M)=0$,
then by theorems of Gromov, Lawson and Stolz,
there exists a metric of positive scalar curvature on~$M$.
If~$M$ is not simply connected,
the problem of the existence of metrics of positive scalar curvature on~$M$
is not yet completely solved.
An overview of related theorems can be found in~\StZwei\
or in~\LM, sections~IV.4--7.
On the other hand,
assume that~$m\ge 3$,
and let~$f\colon M\to\R$ be a function which is negative somewhere on~$M$.
Then there exists a Riemannian metric on~$M$ with scalar curvature~$\kappa=f$
by a theorem of Kazdan and Warner (\KW).

Suppose that there exists a metric~$g$
of positive scalar curvature~$\kappa$ on~$M$.
We call such a metric {\em extremal\/}
if any other metric~$g'$ on~$M$ which is ``larger'' in a suitable sense
has smaller scalar curvature in at least one point of~$M$.
In~\Gromov, section~$5{4\over9}$,
Gromov asked, which manifolds~$M$ possess such an extremal metric,
and how such a metric may look like.
He proposed to specify the word ``larger'' above to mean larger
on two-vectors.
In this case,
we call~$g$ {\em area-extremal.\/}
Gromov also proposed to investigate not only variations of the
metric on~$M$ itself,
but to consider also area-nonincreasing spin maps
of non-vanishing $\Adach$-degree from other Riemannian manifolds to~$M$.
A metric that is area-extremal in this stronger sense will
be called area-extremal {\em in the sense of Gromov.\/}

Up to now,
all known examples of extremal metrics are symmetric:
Llarull showed that the round metrics on spheres are area-extremal
in the sense of Gromov (\LlaZwei).
He also showed that it is not enough to restrict oneself
to metrics which are larger only on three-vectors (\LlaEins).
Min-\penalty1000Oo proved that Hermitian symmetric spaces of compact type
are area-extremal (\MinOo),
and Kramer established extremality of complex and quaternionic
projective spaces (\Kramer).
A general theorem for all compact Riemannian symmetric spaces~$G/K$
with~$\rk G=\rk K$ will be given in~\GS.

In this paper,
we show that K\"ahler metrics of non-negative Ricci curvature on~$M$
are ``\spinc\ area-extremal'' according to \AreaExtremalDefinition.
Here, we compare the scalar curvatures of~$M$
and another Riemannian manifold~$N$
via area-nonincreasing maps~$f\colon N\to M$ of non-vanishing \spinc\ degree,
see~\SpincDegreeFormel.
The following theorem follows from \GromovCorollary\ and \GromovRigidity\
below:

\Theorem\GromovTheorem.
Let~$M$ be a compact, connected K\"ahler manifold
with non-negative Ricci curvature and scalar curvature~$\kappa$,
equipped with its canonical \spinc\ structure.
Let~$N$ be a compact, connected, orientable Riemannian manifold
with scalar curvature~$\kapquer$.
If~$f\colon N\to M$ is an area-nonincreasing spin map
of non-vanishing \spinc\ degree,
then
	$$\kapquer\not>\kappa\circ f\;.$$
If~$g$ is Ricci-positive and~$\kapquer\ge\kappa\circ f$,
then~$\kapquer=\kappa\circ f$.
Moreover, $N$ is then isometric to~$M\times F$,
where~$F$ admits a non-trivial parallel untwisted spinor,
and~$f$ is the projection onto the first factor.

Here~$\kapquer\not>\kappa\circ f$ means
that there exists a~$p\in N$ with~$\kapquer(p)\le\kappa(f(p))$.
For a complex manifold~$M$,
the \spinc\ degree of the identity~$\id_M$ with respect to the natural
\spinc\ structure
is just the Todd genus of~$M$,
which can be zero.
Thus,
in order to conclude that a K\"ahler metric~$g$ of non-negative Ricci
curvature is area-extremal among metrics on~$M$,
we need some extra condition,
cf.\ \AreaExtremalRemark.
On the other hand,
for the estimate in \GromovTheorem\ alone,
we may replace the condition that~$f$ be area-nonincreasing
by a weaker assumption, cf.\ \AreaNonincreasingRemark.

It would be interesting to know the class of complex manifolds
whose Ricci-positive K\"ahler metrics are not only \spinc\ area-extremal,
but also area-extremal in the sense of Gromov.
In this paper,
we establish this property for two series of Hermitian symmetric spaces.
Let
	$$Q^n:=\SO_{n+2}/\SO_n\times\SO_2
	\cong\bigl\{\,[z_0:\ldots:z_{n+1}]
		\bigm|z_0^2+\cdots+z_{n+1}^2=0\,\bigr\}\subset\C P^{n+1}$$
be the complex hyperquadric.
The following estimate generalizes Min-Oo's theorem~7 (\MinOo)
in the special cases of~$\C P^n$ and~$Q^n$:

\Theorem\ZweitesGromovTheorem.
Let~$g$ be a K\"ahler metric on the complex manifold~$M=\C P^n$ or~$M=Q^n$
with non-negative Ricci curvature and scalar curvature~$\kappa$.
Let~$N$ be a compact, connected orientable Riemannian manifold
with scalar curvature~$\kapquer$.
If~$f\colon N\to M$ is an area-nonincreasing spin map
of non-vanishing $\Adach$-degree,
then
	$$\kapquer\not>\kappa\circ f\;.$$
If~$g$ is Ricci-positive and~$\kapquer\ge\kappa\circ f$,
then~$\kapquer=\kappa\circ f$.
Moreover, $N$ is then isometric to~$M\times F$,
where~$F$ admits a non-trivial parallel untwisted spinor,
and~$f$ is the projection onto the first factor.

\Rem
Both \GromovTheorem\ and \ZweitesGromovTheorem\ exhibit
large families of area-extremal metrics
on Ricci-positive K\"ahler manifolds,
cf.\ \forward\ManyMetrics\ and \forward\ManyMetricsZwei\ below.
This is a new phenomenon:
The only examples of area-extremal metrics discovered before
were the standard metrics on certain symmetric spaces,
and thus unique (up to rescaling) on the underlying differentiable manifold.
\endremark

Since all K\"ahler manifolds with positive Ricci curvature admit
a holomorphic (but not necessarily isometric)
embedding into complex projective space,
it is natural to study the scalar curvature~$\kappa$
of smooth projective varieties.
More generally,
one may consider maps~$f$ from a manifold~$M$ into~$\C P^N$.
Recently, B\"ar and Bleecker compared the scalar curvature
of a Riemannian manifold~$M$ with the curvature tensor of another
manifold~$N$
if there exists an area-nonincreasing map~$f\colon M\to N$
of a certain topological type.
For the special case of a complete intersection~$M$
and a map~$f$ into~$\C P^N$ homotopic to the inclusion,
they obtained an explicit upper bound for~$\min\kappa$ (\BB).
In \MainCorollary,
we present a stronger estimate for a certain class of maps into~$\C P^N$,
which in particular implies Conjecture~6.1 in~\BB.
In this estimate, we obtain equality if~$M=\C P^n$ or~$M=Q^n$.

For general~$M$ and~$f$,
\MainCorollary\ is rather coarse
since it does not reflect the topology of~$M$.
If we restrict our attention to complete intersections~$V\subset\C P^N$,
which form a subclass of the class of smooth projective varieties,
we obtain much stronger estimates:

\Theorem\IntroComplIntTheorem.
Let~$V=V^n(a_1,\dots,a_r)$ be a complete intersection
of total degree~$\abs a:=a_1+\dots+a_r$.
Let~$V$ be equipped with a Riemannian metric
such that there exists an area-nonincreasing map~$f\colon V\to\C P^{n+r}$
which is homotopic to the inclusion.
Then
	$$\min_{p\in V}\kappa(p)
	\le\cases
		4n\(n+r+1-\abs a\)
			&\text{if~$\abs a\le n+r$,}\\
		0	&\text{if~$\abs a>n+r$, $n$ is even,
				and~$V$ is spin,}\\
		4n	&\text{if~$\abs a>n+r$,
				and~$V$ is not spin,\quad and}\\
		8n	&\text{if~$\abs a>n+r$, $n$ is odd, and~$V$ is spin.}
	\endcases$$
If~$V$ is connected,
$\abs a\le n+r$
and~$\kappa\ge 4n(n+r+1-\abs a)$,
then~$V$ is a K\"ahler-Einstein manifold
with Einstein constant~$2(n+r+1-\abs a)>0$,
and~$f$ is an isometric, holomorphic embedding.

As the proofs of many of the results cited above,
our arguments are based on a combination of the BLW-formula
with the index theorem,
applied to certain twisted Dirac operators.
The estimates obtained from the BLW-formula become stronger
the simpler the algebraic structure of the twisting curvature
and the smaller its norm is.
In this paper,
we work with the Dirac operators associated to \spinc\ structures,
i.e., we twist with a line bundle
which will be of sufficiently small curvature.
Since the BLW-formula for \spinc\ Dirac operators is particularly simple,
we do not actually need the map~$f$ to be area-nonincreasing
to obtain our estimates.
It suffices that~$\norm{f^*F^L}\le\norm{F^L}\circ f$
with respect to the metrics we compare.
Here~$F^L$ is the curvature of the canonical line bundle of the
\spinc\ structure on~$M$,
and $\norm{F^L}$ is defined in~\OmegaNormDefinition.
We show in \NormLemma\ that indeed~$\norm{f^*F^L}\le\norm{F^L}\circ f$
for any area-nonincreasing~$f$.

We also apply the Atiyah-Singer index theorem
to ensure the existence of a harmonic spinor:
For example,
for K\"ahler manifolds~$M$ of non-negative Ricci curvature 
which is positive somewhere on~$M$,
the Todd genus of~$M$ is~$1$ by Bochner's theorem.
This guarantees the existence of a harmonic spinor
even if the metric is no longer K\"ahler after deformation.

In cases where we do not know the index of a particular \spinc\ Dirac operator
explicitly,
we try to construct a family of \spinc\ structures indexed
by a parameter~$k\in\Z$.
The indices of the corresponding \spinc\ Dirac operators~$D_k$
are given by a polynomial in~$k$,
called the {\em ``Hilbert polynomial'',\/}
which is non-zero under certain mild topological restrictions.
This implies that
we find a harmonic spinor for some~$k$ which is not too large.
For example,
when considering submanifolds of complex projective space,
we construct \spinc\ structures whose canonical line bundle
is a small power~$f^*H^k$ of the pulled-back hyperplane bundle~$H$ of~$\C P^N$.
In the special case of complete intersections,
we can determine explicitly the smallest~$k$
which produces a non-vanishing index,
and obtain the sharper estimate of \IntroComplIntTheorem.

The rest of this paper is organized as follows:
In \EstimateKapitel,
we recall Hitchin's scalar curvature estimate
for \spinc\ Dirac operators with non-vanishing index (\Estimate).
If for some Riemannian manifold~$(M,g)$ one has equality in the estimate
mentioned above,
then~$M$ has special holonomy (\Rigidity).
In \SpinCKapitel,
we define the notion of \spinc\ area-extremality and prove \GromovTheorem.
We introduce the Hilbert polynomial in \HilbertKapitel\
and use it to derive another scalar curvature estimate (\MainTheorem)
and to prove \ZweitesGromovTheorem.
In the last section,
we establish a general scalar curvature estimate for certain submanifolds
of~$\C P^N$ (\MainCorollary).
Using the detailed knowledge of the topology of complete intersections,
we derive the stronger estimate of \IntroComplIntTheorem.
This estimate is then compared with other known results.

Parts of this work were written while the second named author
enjoyed the hospitality and support of the~IHES (Bures-sur-Yvette).
The first named author would like to thank
the Universit\'e de Paris-Sud (Orsay) for its hospitality.
We are grateful to D.~Huybrechts and D.~Kotschick for helpful discussions,
and in particular to C.~B\"ar for bringing
the subject of complete intersections to our attention.

\Section\EstimateKapitel=1
\head 1. Scalar Curvature, Twisted Dirac Operators
and K\"ahler metrics\endhead
In this section,
we recall some facts on \spinc\ structures and parallel spinors
from~\Hitchin, \LM\ and~\Moro.
Let~$M$ be a compact, oriented Riemannian manifold.
Each \spinc\ structure on~$M$ possesses a natural Dirac operator.
Whenever the index of this operator is non-zero,
one obtains an estimate for the scalar curvature of~$M$.
We also investigate the case where equality is obtained
for some metric~$g$ on~$M$:
in this case,
$M$ has to be locally the product of a K\"ahler manifold
with a manifold that admits an untwisted parallel spinor.

We start by recalling the construction of \spinc\ structures
and their associated Dirac operators (\Hitchin, section~1.1, \LM, appendix~D).
Let~$\Spin_2$ be the connected double cover of~$\U_1=\SO_2$,
then the group~$\Spin^c_m$ is defined as
	$$\Spin^c_m:=\Spin_m\times\Spin_2/\{(1,1),(-1,-1)\}\;.$$

\Definition\SpincDefinition
Let~$P_{\SO_m}$ be the frame bundle associated to the tangent bundle~$TM$
of an $m$-dimensional  Riemannian manifold~$M$.
A {\em \spinc\ structure\/} on~$M$ consists
of a principal $\U_1$-bundle~$P_\U$,
a principal $\Spin^c_m$-bundle~$P_\Spinc$,
and a bundle map~$P_\Spinc\to P_{\SO_m}\times P_\U$,
which is equivariant with respect to the natural action of~$\Spin^c_m$
on~$P_{\SO_m}\times P_\U$.
Let~$L$ be the complex line bundle associated to~$P_\U$,
then~$L$ is called the {\em canonical line bundle\/}
associated to~$P_\Spinc$,
and its first Chern class~$c=c_1(L)\in H^2(M,\Z)$
is called the {\em canonical class\/} of~$P_\Spinc$.
\enddefinition

A \spinc\ structure with canonical line bundle~$L$ exists
iff the second Stiefel-Whitney classes of~$M$ and~$L$ coincide,
i.e.\ iff
	$$w_2(M):=w_2(TM)=w_2(L)\in H^2(M,\Z_2)\;.$$
Note that~$w_2(L)$ is just the reduction modulo~$2$ of~$c$.

Let~$S$ be the complex, unitary spinor representation of~$\Spin_m$.
Then~$S$ is irreducible if~$m$ is odd,
and splits as~$S=S^+\oplus S^-$ if~$m$ is even.
The groups~$\Spin_m\times\U_1$ and~$\Spin^c_m$ act on~$S$,
where~$\U_1\subset\C$ acts by complex multiplication.
The {\em complex spinor bundle\/}
associated to a \spinc\ structure~$P_\Spinc$
is the associated Hermitian vector bundle~$\calS$ with fiber~$S$.
The tangent bundle~$TM$ acts on~$\calS$ by Clifford multiplication.

Fix a unitary connection~$\nabla^L$ on the canonical line bundle~$L$,
and let~$F^L\in i\,\Omega^2(M)$ be the imaginary-valued curvature form
of~$\nabla^L$.
Then~$-{1\over2\pi i}F^L$ represents the image~$c_\R$ of~$c$ in~$H^2(M,\R)$,
and each closed two-form on~$M$ which represents~$c_\R$ arises this way.
Together with the Levi-Civita connection,
$\nabla^L$ induces a unitary connection~$\nabla^\calS$ on~$\calS$
which is compatible with Clifford multiplication.
Thus, we can define the {\em \spinc\ Dirac operator\/} associated to~$P_\Spinc$
and~$\nabla^L$ as
	$$D:=\sum_{i=1}^m c_i\nabla_{e_i}^\calS
		\,\colon\;\Gamma(\calS)\longrightarrow\Gamma(\calS)\;.$$
Here, $e_1$, \dots, $e_m$ is a local orthonormal frame of~$M$,
and~$c_i$ denotes Clifford multiplication with~$e_i$.
The following BLW-formula for the square of the Dirac operator~$D$
associated to~$P_\Spinc$ was already known to Schr\"odinger
(\Schro, $\S7$):
	$$D^2=\nabla^*\nabla+{\kappa\over4}
		+{1\over2}\sum_{i<j}F^L(e_i,e_j)\,c_ic_j\;.
	\Formel\SpincBLWFormel$$
Here, $\nabla^*\nabla$ is the Bochner-Laplacian on~$\calS$,
which is a non-negative elliptic differential operator of degree~2,
and~$\kappa$ is the scalar curvature of~$M$.
Recall that on a compact manifold,
$\ker(\nabla^*\nabla)$ is precisely the space of parallel spinors
with respect to~$\nabla^\calS$.

We assume from now on that the dimension of~$M$ is even and equals~$2n$.
Then the complex spinor bundle splits as~$\calS=\calS^+\oplus\calS^-$,
and~$D$ splits as~$D^\pm\colon\Gamma(\calS^\pm)\to\Gamma(\calS^\mp)$.
From the Atiyah-Singer index theorem for \spinc\ Dirac operators
(Theorem~D.15 in~\LM),
we know that
	$$\ind(D)
	:=\dim\ker(D^+)-\dim\ker(D^-)
	=\Bigl(\Adach(M)\,e^{\textstyle{c\over2}}\Bigr)[M]\;.
		\Formel\SpincIndexFormel$$

For the applications in this paper,
we have to measure curvatures of unitary line bundles.
Therefore,
we define a norm on~$\Lambda^2(M)$ by
	$$\norm\alpha_g=\sum_{j=1}^n\abs{\lambda_j}\;,\qquad
		\text{where}\quad
	\alpha=\sum_{j=1}^n\lambda_j\,e^{2j-1}\wedge e^{2j}\in\Lambda^2(M)
		\Formel\OmegaNormDefinition$$
with respect to a suitable $g$-orthonormal frame~$e_1$, \dots, $e_{2n}$
with dual frame~$e^1$, \dots, $e^{2n}$.
From the proof of \NormLemma\ below,
it will become clear that~$\norm\punkt_g$ is indeed a norm.

\Example\KaehlerExample
Let~$(M^{2n},g,J)$ be a K\"ahler manifold,
and let~$\Lambda^{0,*}M$ be the bundle of antiholomorphic differential forms
on~$M$.
Then~$M$ possesses a {\em natural \spinc\ structure\/}~$P_\Spinc M$,
such that the canonical line bundle of~$P_\Spinc M$
is precisely the canonical line bundle~$K=\Lambda^{0,n}M$ of~$M$,
and the canonical class~$c$ of~$P_\Spinc M$ is
the canonical class~$c_1(K)$ of~$M$.
The complex spinor bundle~$\calS$ of this \spinc\ structure
is isomorphic to~$\Lambda^{0,*}M$,
where Clifford multiplication is given
by~$c_{2j-1}+ic_{2j}=-2(\iota_{2j-1}+i\iota_{2j})$
and~$c_{2j-1}-ic_{2j}=\eps^{2j-1}-i\eps^{2j}$.
Here,
$\iota_j$ and~$\eps^j$ denote interior and exterior multiplication with~$e_j$.
The Levi-Civita connection on~$M$
(or equivalently, the holomorphic unitary connection on~$K$)
induces a connection on~$\Lambda^{0,*}M$
such that the associated \spinc\ Dirac operator~$D$ coincides
with the Dolbeault operator~$\bar\delta+\bar\delta^*$.
Thus, \spinc\ Dirac operators can be regarded as natural generalizations
of Dolbeault operators on K\"ahler manifolds.
For K\"ahler manifolds, we have~$\ind(D)=\Td(M)[M]$,
and the constant function~$1$ is always a harmonic
and moreover parallel spinor.
Note also that the curvature of~$K$ is related to the Ricci curvature of~$M$ by
	$$F^K=i\Ric(\punkt,J\punkt)\;,\Formel\RicciFormel$$
so the scalar curvature of~$M$ is given by~$\kappa=\tr(\Ric)\le 2\norm{F^K}_g$,
with equality if the Ricci curvature is non-negative.
\endexample

For a \spinc\ structure with canonical line bundle~$L$
on Riemannian a manifold~$(M,g)$,
there is still a relation between~$\kappa$ and~$\norm{F^L}_g$
if the corresponding \spinc\ Dirac operator has a kernel.
Combining \SpincBLWFormel\ and \SpincIndexFormel\ above,
one obtains

\Proposition\Estimate({\Hitchin, Theorem~1.1}).
Let~$P_\Spinc$ be a \spinc\ structure with canonical class~$c\in H^2(M,\Z)$
on a compact, oriented Riemannian manifold~$(M,g)$,
and let~$c_\R\in H^2(M,\R)$ be represented
by the closed two-form~$-{1\over2\pi}\alpha$.
If the scalar curvature of~$M$ satisfies
	$$\kappa>2\norm\alpha_g$$
everywhere on~$M$,
then~$\bigl(\Adach(M)\,e^{{c\over2}}\bigr)[M]=0$.

\Proof
Choose a connection~$\nabla^L$ on~$L$ with curvature~$F^L=i\alpha$,
and let~$D$ be the Dirac operator associated to~$P_{\Spinc}$ and~$\nabla^L$.
Let~$0\ne\psi\in\Gamma(\calS)$ be a spinor,
and let~$\norm\punkt$ and~$\<\punkt,\punkt\>$ denote the $L^2$-norm
and $L^2$-scalar product, then
\NeueFormel\BLWRechnung
$$\split
	\norm{D\psi}^2
	&=\bigl\<D^2\psi,\psi\bigr\>
	=\biggl<\biggl(\nabla^*\nabla+{\kappa\over4}
			+{1\over2}\sum_{j<k}F^L\(e_j,e_k\)\,c_jc_k\biggr)\psi,
		\psi\biggr\>\\
	&=\bigl\|{\nabla^\calS\psi}\bigr\|^2
		+\Bigl\<{\kappa\over4}\psi,\psi\Bigr\>
		+\,\biggl\<{1\over2}\sum_{j=1}^n
			\alpha\(e_{2j-1},e_{2j}\)\,ic_{2j-1}c_{2j}\,\psi,
			\psi\biggr\>\\
	&\ge\bigl\|{\nabla^\calS\psi}\bigr\|^2
		+\biggl\<\biggl({\kappa\over4}-{\norm\alpha_g\over2}\biggr)
			\psi,\psi\biggr\>
	>0\;,
\endsplit\HierFormel\BLWRechnung$$
Here, we have chosen the frame~$e_1$, \dots, $e_{2n}$ at each point
such that~\OmegaNormDefinition\ holds.
We have also used that~$ic_jc_k$ has eigenvalues~$\pm 1$.
It follows that~$D$ is invertible,
and that its index vanishes.
\qed\enddemo

Now, we want to know for which~$M$ and~$\alpha$ the estimate above
can be sharp.
As already observed by Hitchin,
equality in \Estimate\ implies the existence of a parallel spinor
in~$\Gamma(\calS)$ (\Hitchin).
By~\Moro,
$M$ has special holonomy in this case.
If moreover, $\alpha^n\ne 0$ at some point of~$M$,
then one can conclude that~$M$ is K\"ahler:

\Theorem\Rigidity.
Let~$(M^{2n},g)$ be a compact, connected, orientable Riemannian manifold,
let~$P_\Spinc$ be a \spinc\ structure on~$M$ with canonical class~$c$,
and let~$c_\R\in H^2(M,\R)$ be represented
by the closed two-form~$-{1\over2\pi}\alpha$.
Assume that~$\bigl(\Adach(M)\,e^{{c\over2}}\bigr)[M]\ne 0$.
If the scalar curvature of~$M$ satisfies
	$$\kappa\ge2\norm\alpha_g$$
everywhere on~$M$,
then~$\kappa=2\norm\alpha_g$,
and the universal cover of~$M$ is the product of a K\"ahler manifold
with a manifold which admits an untwisted parallel spinor.
\parr
Assume moreover that~$\alpha^n\ne 0$ in at least one point of~$M$.
In this case,
there exists a complex structure~$J$
such that~$(M,g,J)$ is a K\"ahler manifold
with natural \spinc\ structure~$P_{\Spinc}$,
and~$i\alpha$ is the curvature of the associated line bundle.
The Ricci curvature of~$M$ is non-negative
and positive somewhere on~$M$,
and the Todd genus of~$M$ is~$1$.

\Remark\ConverseRemark
Conversely,
assume that the Ricci curvature of a connected K\"ahler manifold~$M$ is
non-negative and positive somewhere on~$M$.
Then~$\kappa=2\norm\rho_g$,
where~$\rho$ is the curvature of the canonical line bundle,
and~$\Td(M)[M]=1$ by Bochner's theorem
(\LM, Corollary~IV.11.12).
\endremark

\demo{Proof\/ \rm of \Rigidity}
Let~$\nabla^L$ and~$D$ be as in the proof of \Estimate\ above.
Since~$\ind(D)\ne0$,
there exists a harmonic spinor~$0\ne\psi\in\ker(D)\subset\Gamma(\calS)$.
From \BLWRechnung,
we conclude that
\NeueFormel\BLWRigidity
$$\split
	0=\norm{D\psi}^2
	&=\bigl\|{\nabla^\calS\psi}\bigr\|^2
		+\Bigl\<{\kappa\over4}\psi,\psi\Bigr\>
		+\,\biggl\<{1\over2}\sum_{j=1}^n
			\alpha\(e_{2j-1},e_{2j}\)\,ic_{2j-1}c_{2j}\,\psi
			,\psi\biggr\>\\
	&=\bigl\|{\nabla^\calS\psi}\bigr\|^2
		+\biggl\<\biggl({\kappa\over4}-{\norm\alpha_g\over2}\biggr)
			\psi,\psi\biggr\>\;.
\endsplit\HierFormel\BLWRigidity$$
In particular,
$\psi$ is parallel with respect to~$\nabla^\calS$,
and~$i\lambda_jc_{2j-1}c_{2j}\psi=-\abs{\lambda_j}\psi$,
where at each point,
the frame~$e_1$, \dots, $e_{2n}$ is chosen such that~\OmegaNormDefinition\
holds.

The proof proceeds from now on as in~\Moro, chapter~3.
We will nevertheless give a detailed argument here,
because we will need the explicit construction of the complex structure~$J$
again in the proof of \GromovRigidity\ below.
Suppose that~$\alpha^n\ne0$ at some point~$p\in M$.
This implies that each of the~$\lambda_j$ is non-zero,
and moreover negative
if the frame~$e_1$, \dots, $e_{2n}$ is chosen accordingly.
Thus at~$p$,
	$$ic_{2j-1}c_{2j}\psi=\psi
		\quad\Longrightarrow\quad
	\(c_{2j-1}+ic_{2j}\)\psi=0\;,
		\Formel\PureRelation$$
i.e., $\psi$ is a {\em pure spinor,\/}
cf.~\LM, Definition~IV.9.3.
Now, \Rigidity\ is proved in the same way as Proposition~IV.9.8 in~\LM:
since~$\psi$ is parallel,
\PureRelation\ still holds if the frame~$e_1$, \dots, $e_{2n}$ at~$p$
is parallelly translated to any point of~$M$ along any curve.
Here, we need in particular that~$M$ is connected.
Hence, there is a unique parallel almost complex structure~$J$ on~$M$ such that
	$$\(c(v)+i\,c(Jv)\)\psi=0
		\Formel\KaehlerRelation$$
for all~$v\in TM$,
so~$g$ is a K\"ahler metric.
The definition of~$J$ may be rephrased by saying that
Clifford multiplication with antiholomorphic vectors vanishes on~$\psi$.

In order to prove that~$P_\Spinc$ is the natural \spinc\ structure
associated to~$J$,
we show that~$\calS$ is isomorphic as a Dirac bundle
to the bundle~$\Lambda^{0,*}M$ of antiholomorphic forms on~$M$.
We claim that there is a natural isomorphism which maps the form~$1$ to~$\psi$,
and more generally,
	$$\(e^{2j_1-1}-ie^{2j_1}\)\wedge\dots
		\wedge\(e^{2j_k-1}-ie^{2j_k}\)
	\quad\longmapsto\quad
	\(c_{2j_1-1}-ic_{2j_1}\)\cdots\(c_{2j_k-1}-ic_{2j_k}\)\,\psi
		\Formel\IsomorphismFormel$$
for all~$k\le n$ and~$1\le j_1<\dots<j_k\le n$.
Using~\KaehlerRelation,
one easily checks that this defines an isomorphism of Clifford modules
which is compatible with the Levi-Civita connection~$\nabla$
on~$\Lambda^{0,*}M$ and with~$\nabla^\calS$ on~$\calS$.

Once we have identified the canonical line bundle~$L$
of our given \spinc\ structure with the canonical line bundle of~$(M,g,J)$,
we know that $\alpha$ is the Ricci form of~$M$,
and the Ricci curvature is given by~\RicciFormel.
On the other hand, from
	$$\trace(\Ric)=\kappa=2\norm\alpha_g=2\norm{F^L}_g$$
we can conclude that the Ricci curvature is non-negative,
and positive wherever~$\alpha^n\ne 0$.
By Bochner's theorem, the Todd genus of~$M$ is~$1$.
\qed\enddemo

One might hope that the conditions of \Rigidity\ determine
the K\"ahler metric~$g$ uniquely.
Unfortunately, this is not true in general:

\Example\FlagExample
Let~$G$ be a semi-simple Lie group of rank~$k$ with maximal torus~$T$.
Then~$G/T$ carries as many different homogeneous complex structures
as the Weyl group of~$G$ has elements.
Let us fix a homogeneous complex structure~$J$ and a positive constant~$C$.
Then there exists a $(k-1)$-parameter family of homogeneous K\"ahler metrics
with constant scalar curvature~$\kappa=C$.
All these metrics have the same Ricci curvature~$\Ric$ (\Besse, chapter~8),
which is positive.
Thus,
the Ricci form~$\alpha=\Ric(\punkt,J\punkt)$ depends only on~$J$ and~$\kappa$,
but not on the metric~$g$.
This means that~$g$ cannot be completely determined by~$\alpha$ and~$P_\Spinc$.
\endexample

In order to get a better control of the metric~$g$,
one has to introduce more constraints.
For example,
one could restrict oneself to metrics which are in a suitable sense
``larger'' than a fixed background metric.
This is precisely what we do in \GromovRigidity\ in the next section
in a slightly more general setting.

\Section\SpinCKapitel=2
\head 2. \SPinc\ Area-Extremality of Ricci-positive K\"ahler Metrics\endhead
We apply the results of the previous section to smooth, area-nonincreasing maps
from some Riemannian manifold~$(N,\gquer)$ to a K\"ahler manifold~$(M,g,J)$
of positive Ricci curvature.
We show that a K\"ahler metric of positive Ricci curvature
is ``area-extremal'' in a sense closely related to
Gromov's definition of area-extremality in~\Gromov.

We recall that
a map~$f\colon N\to M$ is called a {\em spin map\/} iff
	$$w_2(N)=f^*w_2(M)\;.$$
Let~$P_\Spinc$ be \spinc\ structure on~$M$ with canonical class~$c$.
Then we define the {\em \spinc\ degree\/} of~$f$
to be
	$$\bigl(\Adach(N)\,f^*e^{c\over2}\bigr)[N]\;.
		\Formel\SpincDegreeFormel$$
Clearly, this is an integer if~$f$ is a spin map.
If~$f$ is a covering map,
then~$\Adach(N)=f^*\Adach(M)$,
and the \spinc\ index
is precisely~$\deg f\cdot\bigl(\Adach(M)\,e^{c\over2}\bigr)[M]$.
Thus, the \spinc\ degree generalizes the degree of a covering map
whenever~$\bigl(\Adach(M)\,e^{c\over2}\bigr)[M]=1$.
Let~$g$ and~$\gquer$ be Riemannian metrics on~$M$ and~$N$.
Then~$f$ is called {\em area-nonincreasing\/}
iff
	$$\abs{f_*(v\dach w)}_g\le\abs{v\dach w}_{\gquer}$$
for all decomposable two-vectors~$v\dach w\in\Lambda^2TN$.
If~$f$ and~$h$ are two functions on~$M$,
we write
	$$f\not>h$$
iff there exists a~$p\in M$ with~$f(p)\le h(p)$,
i.e., iff~$f$ is {\em not everywhere\/} greater than~$h$.

\Definition\AreaExtremalDefinition
Let~$(M,g)$ be compact, connected, oriented Riemannian manifold
with a \spinc\ structure with canonical class~$c$,
and let~$\kappa$ be the scalar curvature of~$M$.
Then the metric~$g$ is called {\em \spinc\ area-extremal\/}
iff
	$$\kapquer\not>\kappa\circ f$$
for all compact, connected, oriented Riemannian manifolds~$(N,\gquer)$
with scalar curvature~$\kapquer$
and all smooth, area-nonincreasing spin maps~$f\colon N\to M$
of non-vanishing \spinc\ degree.
\enddefinition

\Remark\AreaExtremalRemark
This notion of \spinc\ area-extremality is closely related
to area-extremality in the sense of Gromov (\Gromov, section~$5{4\over9}$).
Recall that the {\em $\Adach$-degree\/} of a map~$f\colon N\to M$
is defined as~$\(\Adach(N)\,f^*v\)[N]$,
where~$v$ is the canonical generator of~$H^m(M,\Z)$
for a compact, connected, oriented $m$-di\-men\-si\-o\-nal manifold~$M$.
Gromov calls a metric~$g$ {\em area-extremal\/}
if the scalar curvatures~$\kapquer$ of~$N$
is smaller or equal to~$\kappa\circ f$ in at least one point of~$N$
for all area-nonincreasing spin maps of non-vanishing $\Adach$-degree.
Because the $\Adach$-degree generalizes the degree of a map between
compact, connected, oriented manifolds of the same dimension,
the identity of~$M$ always has $\Adach$-degree~$1$,
and one can compare any metric on~$M$ with an area-extremal metric
using~$f=\id_M$.

On the other hand,
if we want to use \spinc\ area-extremality
to compare different metrics on~$M$,
we need~$M$ to be connected,
and~$\id_M$ to have non-zero \spinc\ degree.
By definition,
this means that~$\bigl(\Adach(M)\,e^{c\over2}\bigr)[M]\ne 0$.
If we regard the natural \spinc\ structure of a complex manifold,
then the \spinc\ degree of~$\id_M$ is just the Todd genus of~$M$,
which can be zero.
In particular, \spinc\ area-extremality does not automatically imply
area-extremality.
However,
if~$M$ carries a K\"ahler metric~$g$ of non-negative Ricci curvature
which is positive somewhere,
the \spinc\ degree of~$\id_M$ is~$1$ by Bochner's theorem
(see \ConverseRemark),
so~$g$ is indeed area-extremal.
\endremark

From \Estimate,
we get immediately a large class of \spinc\ area-extremal metrics:

\Corollary\GromovCorollary.
Let~$(M,g,J)$ be a compact, connected K\"ahler manifold
of non-negative Ricci curvature.
Then the metric~$g$ on~$M$ is \spinc\ area-extremal.
If the Ricci curvature is positive somewhere on~$M$,
then~$g$ is also area-extremal.

\Proof
Let~$K\to M$ be the canonical line bundle with curvature~$F^K=i\alpha$,
and let~$f^*K$ be the pull-back of~$K$ to~$N$,
equipped with the pull-back connection.
Because~$f$ is a spin map,
$N$ carries a \spinc\ structure with canonical line bundle~$f^*K$.
The index of the associated Dirac operator is precisely the \spinc\ degree
of the map~$f$.
We prove in \NormLemma\ below
that~$\|f^*\alpha\|_\gquer\le\norm\alpha_g\circ f=\kappa\circ f$
if~$f$ is area-nonincreasing.
Now, the first claim follows by applying \Estimate\ to~$N$ and~$\alpha$.
For the second statement,
one uses Bochner's theorem as in the remark above.
\qed\enddemo

\Remark\ManyMetrics
Let~$(M,g,J)$ be a compact, connected K\"ahler manifold
with non-negative Ricci curvature which is positive somewhere.
Then not only is~$(M,g)$ area-extremal,
but the same holds for any K\"ahler metric~$g'$ on~$M$ which is $C^2$-close
to~$g$,
such that~$\supp(g'-g)$ is contained in the region~$M^+\subset M$
where the Ricci curvature of~$(M,g)$ is positive.
Such metrics can be constructed
by taking a $C^4$-small function~$h\colon M\to\R$ with~$\supp(h)\subset M^+$:
Let~$\omega$ be the K\"ahler form of~$(M,g,J)$,
then~$\omega+i\,\del\delquer h$ is the K\"ahler form of a metric~$g'$
as above.
This shows that on all compact, connected K\"ahler manifolds
with non-negative Ricci curvature which is positive somewhere,
there exists an infinite-dimensional family of area-extremal metrics.
\endremark

\Remark\AreaNonincreasingRemark
Clearly, the statement of the corollary remains correct if
in \AreaExtremalDefinition,
the condition that~$f$ is area-nonincreasing
is replaced by the much weaker requirement
that~$\norm{\smash{f^*\alpha}}_\gquer\le\norm\alpha_g\circ f$.
The full power of the stronger condition in \AreaExtremalDefinition\
will only become visible in \GromovRigidity\ below,
where we investigate the case that~$\kapquer\ge\kappa\circ f$.
\endremark

Let us check that
indeed~$\norm{\smash{f^*\alpha}}_\gquer\le\norm\alpha_g\circ f$
if~$f$ is area-nonincreasing:

\Lemma\NormLemma.
Let~$(M^{2n},g)$ be a Riemannian manifold,
and let~$\norm\punkt_g$ be the norm on~$\Lambda^2M$
introduced in~\OmegaNormDefinition.
If~$(N,\gquer)$ is another Riemannian manifold
and~$f\colon N\to M$ is a smooth, area-nonincreasing map,
then
	$$\norm{f^*\alpha}_{\gquer}
	\le\norm\alpha_g\circ f$$
for all alternating two-forms~$\alpha$.
Assume that~$n\ge 2$ and that~$\alpha^n\ne 0$ everywhere on~$M$,
then equality implies that~$f$ is a Riemannian submersion.

\Proof
We show first that~$\norm\alpha_g$
can equivalently be defined as
	$$\norm\alpha_g
	=\max\biggl\{\,\sum_{j=1}^n\alpha\(e_{2j-1},e_{2j}\)
		\biggm|\text{$e_1, \dots, e_{2n}$ is an orthonormal
			frame with respect to~$g$}\,\biggr\}\;.
	\Formel\AlternativeDefinition$$
Let~$\alpha=a_1e^1\dach e^2+\dots+a_ne^{2n-1}\dach e^{2n}$
be an alternating form on~$\R^{2n}$,
equipped with the standard Euclidean metric~$g$.
Let~$J$ be the matrix of the standard complex structure on~$\R^{2n}$,
and let~$A$ be the skew-symmetric matrix that represents~$\alpha$.
Then~\AlternativeDefinition\ can clearly be rewritten as
	$${1\over2}\,\max\bigl\{\,\trace\(BAB^{-1}J\)
		\bigm|B\in\O(2n)\,\bigr\}\;.$$

Writing~$B_t=e^{tX}B$ for a skew-symmetric matrix~$X$,
we see that~$B$ is critical point
of the functional~$B\mapsto\trace\(BAB^{-1}J\)$ iff
	$$0={\del\over\del t}\Big|_{t=0}\trace\(B_tAB_t^{-1}J\)
	=\trace\(XBAB^{-1}J\)-\trace\(BAB^{-1}XJ\)
	=2\trace\(XBAB^{-1}J\)$$
for all~$X\in\so(2n)$.
This happens iff~$BAB^{-1}J$ is a symmetric matrix,
i.e.\ iff~$J$ and~$BAB^{-1}$ commute.
In this case, the eigenvalues of~$BAB^{-1}J$
are precisely~$\pm a_1$, \dots, $\pm a_n$,
each with multiplicity~$2$.
The maximum clearly occurs if all these eigenvalues are positive.
For such a~$B$,
${1\over2}\trace\(BAB^{-1}J\)$
is indeed the norm defined in~\OmegaNormDefinition.
This proves~\AlternativeDefinition.
Note that from~\AlternativeDefinition,
one easily concludes that~$\norm\punkt_g$ is sub-additive,
so it is indeed a norm.

Let~$\alpha=\alpha_1+\dots+\alpha_n$
with~$\alpha_i=a_ie^{2i-1}\dach e^{2i}$.
Then~$\norm\alpha_g=\norm{\smash{\alpha_1}}_g+\dots+\norm{\smash{\alpha_n}}_g$.
Suppose that~$f\colon\R^l\to\R^{2n}$ is a linear area-nonincreasing map
with respect to the standard Euclidean metrics~$\gquer$ on~$\R^l$
and~$g$ on~$\R^{2n}$.
Then~$\norm{f^*(e^j\dach e^k)}_{\gquer}\le\norm{e^j\dach e^k}_g$.
Thus, we have
	$$\norm{f^*\alpha}_{\gquer}
	\le\norm{f^*\alpha_1}_{\gquer}+\dots+\norm{f^*\alpha_1}_{\gquer}
	\le\norm{\alpha_1}_g+\dots+\norm{\alpha_1}_g
	=\norm\alpha_g\;,\Formel\NormEstimate$$
which proves the inequality in the lemma.

Let us assume that equality holds,
and that~$n\ge 2$ and~$\alpha^n\ne 0$.
Splitting~$\alpha=\alpha_1+\dots+\alpha_n$ as above,
we find that by~\NormEstimate,
there are orthonormal frames~$\equer_1$, \dots, $\equer_l$ of~$\R^l$
and~$e_1$, \dots, $e_{2n}$ of~$\R^{2n}$
such that
	$$f^*\alpha=\sum_{j=1}^n\lambda_j\,\equer^{2j-1}\wedge\equer^{2j}
		\qquad\text{and}\qquad
	\alpha=\sum_{j=1}^n\lambda_j\,e^{2j-1}\wedge e^{2j}$$
with~$\lambda_1$, \dots, $\lambda_n>0$,
and moreover,
$f_*(\equer_{2j-1}\wedge\equer_{2j})=e_{2j-1}\wedge e_{2j}$
for~$1\le j\le n$.
We may choose the frames above such that~$f_*\equer_{2j-1}=\mu_je_{2j-1}$
and~$f_*\equer_{2j}=\mu_j^{-1}e_{2j}$, with~$\mu_j\ge 1$.
Because~$f$ is area-nonincreasing,
	$$\mu_j\mu_k=\abs{f_*(\equer_{2j-1}\dach\equer_{2k-1})}_g\le1$$
for~$1\le j<k\le n$.
If~$n\ge 2$, this implies clearly that~$\mu_1=\dots=\mu_n=1$.
Using once more that~$f$ is area-nonincreasing,
one proves that~$f_*\equer_k=0$ for~$k>2n$.
This implies the rigidity statement in \NormLemma.
\qed\enddemo

\Rem
On odd-dimensional manifolds,
the estimate of the lemma holds unchanged.
However,
in the case of equality,
one gets a weaker statement
because one cannot control~$f_*$ on~$\ker(f^*\alpha)$.
\endremark

The rigidity statement of \NormLemma\
can be used to investigate the case
that we have the inequality~$\kapquer\ge\kappa\circ f$
in \AreaExtremalDefinition:

\Theorem\GromovRigidity.
Let~$(M,g)$ be a compact, connected K\"ahler manifold
of positive Ricci curvature and complex dimension~$n\ge 2$,
let~$(N,\gquer)$ be another compact, connected, oriented Riemannian manifold,
and let~$f\colon N\to M$ be a smooth spin map of non-zero \spinc\ degree.
Suppose that~$f$ is area-nonincreasing
and that~$\kapquer\ge\kappa\circ f$.
Then~$N$ is a Riemannian product~$N=M\times F$,
and~$f$ is the projection onto the first factor.
The manifold~$F$ carries a parallel untwisted spinor and is in particular
Ricci flat.

\Proof
Let~$K$ be the canonical line bundle of~$M$
with curvature~$F^K=i\alpha$.
As in the proof of \GromovCorollary,
we construct a \spinc\ structure~$P_\Spinc$ on~$N$
with canonical line bundle~$f^*K$ which has curvature~$f^*F^K$.
If~$f$ is area-nonincreasing and~$\kapquer\ge\kappa\circ f$,
then
	$$2\norm{\smash{f^*\alpha}}_\gquer=\kapquer
		=\kappa\circ f=2\norm\alpha_g\circ f$$
as in the proof of \Rigidity.
Because~$M$ is Ricci-positive,
we have~$\alpha^n\ne 0$,
and~$f$ is a Riemannian submersion by \NormLemma.
Moreover,
the complex spinor bundle associated to~$P_\Spinc\to N$
has a non-trivial parallel spinor.
Then  by~\Moro,
the universal cover~$\pi\colon \Nsl\to N$
splits as a product of a K\"ahler manifold~$\Nsl_1$
with a manifold~$\Nsl_2$
which admits a parallel untwisted spinor and contains no
de Rham factor that is K\"ahler.
We will prove first that~$f\circ\pi$ factors
over the projection~$\pisl_1\colon\Nsl\to \Nsl_1$.
In a second step, we show that~$\Nsl_1$ is isometric to~$M\times F_1$,
such that~$f$ is the projection onto the first factor.

To prove that~$f\circ \pi$ factors over~$\pisl_1$,
let~$\psi$ be the parallel spinor on~$\Nsl$.
We consider the subspace
	$$E:=\bigl\{\,\vquer\in T\Nsl\bigm|
		\text{there is a vector~$\wquer\in T\Nsl$
			such that } c_{\vquer+i\wquer}\psi=0\,\bigr\}$$
as in~\Moro.
Because~$\psi$ is parallel,
this is a parallel distribution in~$T\Nsl$,
and it carries a parallel complex structure~$\Jquer$
such that~$c_{\vquer+i\Jquer\vquer}\psi=0$ for all~$\vquer\in E$
as in~\KaehlerRelation.
Then~$E$ is the tangent distribution to the K\"ahler factor~$\Nsl_1$ of~$\Nsl$.
Let~$\equer_1$, \dots, $\equer_{\nquer}$ be an orthonormal
frame of~$T_p\Nsl$ such that
	$$\pi^*f^*\alpha
	=\sum_{k=1}^n\lambda_k\equer^{2k-1}\wedge\equer^{2k}\;,$$
where~$\lambda_1$, \dots, $\lambda_n$ are negative
as in the proof of \Rigidity.
Thus, $\equer_{2n+1}$, \dots, $\equer_{\nquer}\in\ker(f_*)$.
Because of the analog of equation~\BLWRigidity\ for~$\psi$ on~$\Nsl$,
we conclude that~$\equer_1$, \dots, $\equer_{2n}\in E$,
and that~$\Jquer\equer_{2k-1}=\equer_{2k}$ for~$1\le k\le n$.
This implies already
that~$f\circ\pi$ factors over~$\pisl_1$.

Let~$i_1\colon \Nsl_1\to N$ be an embedding of~$\Nsl_1$ as a factor of~$N$,
and write~$\fquer:=f\circ\pi\circ i_1\colon \Nsl_1\to M$.
Now, $e_1:=\fquer_*\equer_1$, \dots, $e_{2n}:=\fquer_*\equer_{2n}$
form an orthonormal frame of~$T_{f(p)}M$ with
	$$\alpha=\sum_{k=1}^n\lambda_k e^{2k-1}\wedge e^{2k}\;.$$
Because~$\alpha$ is the Ricci form of~$M$,
it follows that~$Je_{2k-1}=e_{2k}$.
To be more precise,
write~$\alpha$ and~$J$ as matrices with respect to the given frame.
Because~$\Ric_M=\alpha(J\punkt,\punkt)$ is symmetric,
these matrices commute.
Finally,
because~$\Ric_M$ is positive and all the~$\lambda_j$ are negative,
$J$ acts as indicated.

In particular, $\fquer$ is a holomorphic map.
Then the fibers~$F_q:=\fquer^{-1}(q)$
are complex submanifolds of~$\Nsl_1$ for all~$q\in M$,
so their mean curvature vector vanishes.
The O'Neill formulas for the Ricci curvature
(\Besse, Proposition~9.36)
imply for the horizontal lift~$X\in T\Nsl_1$ of a vector in~$TM$:
	$$\Ric_{\Nsl_1}(X,X)
	=\Ric_M(\fquer_*X,\fquer_*X)
		-2\sum_{j=1}^{2n}\norm{A(X,\equer_j)}^2
		-\sum_{j=2n+1}^{2n_1}\norm{T(\equer_j,X)}^2\;,
	\Formel\ONeillFormel$$
where~$A(X,Y)={1\over2}[X,Y]^\perp$ is the curvature of the fibre
bundle~$\fquer\colon \Nsl_1\to M$,
and~$T(U,X)=\(\nabla_U^{\Nsl_1}X\)^\perp$ is the shape tensor
of the fibres.
Here, $Y$ is the horizontal lift of another vector in~$TM$,
$U$ is any vertical vector,
$(\punkt)^\perp$ denotes projection onto the space of vertical vectors,
and~$2n_1$ is the dimension of~$\Nsl_1$.

On the other hand,
the Ricci curvatures of~$\Nsl_1$ and~$M$ are related by the equation
	$$\Ric_{\Nsl_1}=-i\(\fquer^*F^K\)\(\Jquer\punkt,\punkt\)
	=\(\fquer^*\alpha\)(\Jquer\punkt,\punkt)
	=\fquer^*\(\alpha(J\punkt,\punkt)\)
	=\fquer^*\Ric_M\;.$$
By \ONeillFormel,
this implies that both tensors~$A$ and~$T$ vanish,
so~$\fquer\colon \Nsl_1\to M$ is locally isometric to a trivial bundle
(\Besse, 9.26).
By a theorem of Kobayashi (\Besse, Theorem~11.26),
$M$ is simply connected,
so this bundle is also globally isometric to a product,
and~$\Nsl_1=M\times F_1$,
where~$\fquer$ is the projection onto the first factor.
In the same way,
$N$ also splits as~$M\times F$ with~$F=F_1\times \Nsl_2$,
such that~$f$ is the projection onto the first factor.
Finally, $F_1$ is a Ricci flat K\"ahler manifold,
so it carries a parallel untwisted spinor.
Since~$\Nsl_2$ also has a parallel untwisted spinor,
the same holds for~$F$.
\qed\enddemo

\Section\HilbertKapitel=3
\head 3. The Hilbert Polynomial
and Area-Extremality \`a la Gromov\endhead
Suppose that~$M^{2n}$ admits a \spinc\ structure
with canonical line bundle~$L$.
Then~$M$ also admits \spinc\ structures that have certain tensor powers~$L^k$
of~$L$ as their canonical line bundle.
The index of the associated Dirac operators~$D_k$ is a polynomial~$P_L(k)$
in~$k$,
which we will call the ``Hilbert polynomial''
in analogy with~\LM, section~IV.11.
If this polynomial is not identically zero,
then it has at most~$n$ zeros.
In this case,
there exists a  small~$k$ such that~$D_k$ has non-vanishing index,
and we can apply the results of the previous chapter
to obtain an estimate of the scalar curvature of~$M$.

For general~$M$,
we do not obtain the best estimate possible,
because the topology of~$M$ is used only superficially.
For complex projective spaces and complex hyperquadrics however,
the estimate is sharp.
For these manifolds,
we recover estimates and rigidity statements as in \SpinCKapitel,
but with the \spinc\ degree replaced by the $\Adach$-degree.
In particular, we prove \ZweitesGromovTheorem.

\Theorem\MainTheorem.
Let~$(M^{2n},g)$ be a compact, connected, orientable Riemannian manifold
of real dimension~$2n$.
Suppose that~$w_2(M)$ is a multiple
of the reduction modulo~$2$ of some~$c\in H^2(M,\Z)$ with~$c^n[M]\ne0$,
and let~$c_\R\in H^2(M,\R)$ be represented
by the closed two-form~$-{1\over2\pi}\alpha$.
Then
	$$\kappa
	\not>\cases
		2n\,\norm\alpha_g	&\text{if~$w_2(M)=nc\bmod2$,
						\quad and}\\
		2(n+1)\,\norm\alpha_g	&\text{if~$w_2(M)=(n+1)c\bmod2$.}
	\endcases$$
Assume that~$\kappa\ge2n\,\norm\alpha_g$ in the first case
or~$\kappa\ge2(n+1)\,\norm\alpha_g$ in the second case.
Then equality holds,
and~$M$ is K\"ahler and biholomorphic
to the complex quadric~$Q^n$ in the first case,
and to~$\C P^n$ in the second.

Recall that~$Q^n\cong\SO_{n+2}/\SO_n\times\SO_2\subset\C P^{n+1}$
was defined by the equation~$z_0^2+\dots+z_{n+1}^2=0$.

\demo{Proof\/ \rm of \MainTheorem}
If~$w_2(M)$ is a multiple of~$w_2(L)=c\bmod2\in H^2(M,\Z_2)$,
we may construct \spinc\ structures~$P_\Spinc^k$ on~$M$
with canonical line bundle~$L^k$ for all~$k\in\Z$
such that~$w_2(M)=kc\bmod 2$.
Let~$D_k$ be the Dirac operator associated to~$P_\Spinc^k$
and~$\nabla^{L^k}$.
By \Estimate,
we know that~$\ind(D_k)=0$
if the scalar curvature~$\kappa$
is everywhere larger than~$2\,\bigl\|F^{L^k}\bigr\|_g=2k\norm\alpha_g$.

On the other hand,
we have assumed that~$c^n[M]\ne 0$.
Then the {\em ``Hilbert polynomial''\/}
	$$P_L(k):=\ind\(D_k\)
	=\Bigl(\Adach(M)\,e^{\textstyle{kc\over2}}\Bigr)[M]$$
is a polynomial of degree~$n$ in~$k$.
It has the non-vanishing leading term~${k^n\over 2^nn!}c^n[M]$,
because~$\Adach(M)$ always starts with~$1$ in degree zero.
In particular, $P_L(k)$ has at most $n$ different zeros.

We distinguish two cases:
if~$w_2(M)=nc\bmod2$,
then for the $n+1$ different values~$n$, $n-2$, \dots, $-n$ of~$k$
there exists a \spinc\ structure~$P_\Spinc^k$,
and the operator~$D_k$ is well-defined.
In particular, for one~$k_0$ with~$\abs{k_0}\le n$,
we have~$\ind(D_{k_0})\ne 0$,
so~$\kappa\not>2\abs{k_0}\norm\alpha_g\le 2n\norm\alpha_g$ by~\Estimate.

For~$w_2(M)=(n+1)\,f^*w_2(H)$, the proof is completely analogous,
but we have to choose~$k_0$ among the~$n+2$
different values~$n+1$, $n-1$, \dots, $-n-1$,
which accounts for the slightly weaker estimate.

We will now study the cases where ``$\ge$'' holds
for the estimate in the theorem.
If we have~$w_2(M)=nc\bmod2$ and~$\kappa\ge 2n\norm\alpha_g$,
then we know not only that~$\kappa=2n\norm\alpha_g$
and that~$M$ is K\"ahler by \Rigidity,
but also that the canonical line bundle~$K$ of~$M$
coincides with~$L^n$.
From a result of Kobayashi and Ochiai,
it follows that~$M$ is biholomorphic to~$Q^n$ (\KO).
Similarly, if~$w_2(M)=(n+1)c\bmod 2$
and~$\kappa\ge 2(n+1)\norm\alpha_g$,
then~$K=L^{n+1}$,
and~$M$ is biholomorphic to~$\C P^n$ (\KO).
\qed\enddemo

Assume that~$M=\C P^n$.
In this case,
we immediately recover the Hilbert-polynomial of~$\C P^n$
from the proof of \MainTheorem:

\Remark\CPnRemark
Let~$M=\C P^n$ be equipped with a K\"ahler metric of positive Ricci curvature.
Then~$D_k$ is invertible for~$k=n-1$, $n-3$, \dots, $1-n$
by \Estimate.
This implies that the zeros of the Hilbert polynomial
	$$P_n(k)
	:=\ind\(D_k^{\C P^n}\)
	=\biggl(\Adach(\C P^n)\,e^{\textstyle{kc\over2}}\biggr)
		[\C P^n]$$
of~$H\to\C P^N$ are precisely~$k=n-1$, $n-3$, \dots, $1-n$
(cf.~\LM, chapter~IV.11).
Since the Todd genus~$P_n(n+1)$ of~$\C P^n$ equals~$1$,
this completely determines~$P_n(k)$.
In particular, its leading term is positive.
We conclude that if~$k\equiv n+1\pmod2$, then
	$$P_n(k)=0\quad\text{if~$\abs{k}\le n-1$,\qquad and}\qquad
	P_n(k)>0\quad\text{if~$k> n-1$.}
	\Formel\CPnIndexEstimate$$
Moreover, since~$\Adach(\C P^n)\in H^{4*}(\C P^n)$,
the polynomial~$P_n(k)$ is even if~$n$ is even, and odd if~$n$ is odd.
\endremark

The theorem above admits a reformulation for
spin maps of non-vanishing $\Adach$-degree.
Recall that the $\Adach$-degree of a map~$f\colon N\to M$
between compact oriented manifolds was defined by
	$$\deg_{\Adach}f:=\(\Adach(N)\,f^*v_M\)[N]\;,$$
where~$v_M$ is the canonical generator of~$H^{\dim M}(M,\Z)$.

\Corollary\AdachDegreeCorollary.
Let~$M^{2n}$ be a compact, connected, oriented manifold,
let~$(N,\gquer)$ be a compact, connected, oriented Riemannian manifold,
and let~$f\colon N\to M$ be a smooth map of non-vanishing $\Adach$-degree.
Suppose that~$w_2(N)$ is a multiple
of the reduction of~$f^*c$ modulo~2
for some~$c\in H^2(M,\Z)$ with~$c^n[M]\ne0$.
Let~$f^*c_\R\in H^2(N,\R)$ be represented
by a closed two-form~$-{1\over2\pi}\alquer$.
Then for the scalar curvature~$\kapquer$ of~$N$,
we have
	$$\kapquer
	\not>\cases
		2n\,\norm{\alquer}_\gquer	&\text{if~$w_2(M)=nc\bmod2$,
							\quad and}\\
		2(n+1)\,\norm{\alquer}_\gquer	&\text{if~$w_2(M)
							=(n+1)c\bmod2$.}
	\endcases$$

\Proof
As in the proof of \MainTheorem,
we can construct \spinc\ structures $P_\Spinc^k$
whenever~$w_2(N)=kf^*c\bmod2$.
The index of the corresponding Dirac operator~$D_k$ is then given by
	$$P_L(k)=\biggl(\Adach(N)\,e^{k\,f^*c\over2}\biggr)[N]\;.$$
Because~$(f^*c)^\nu=f^*(c^\nu)=0$ for~$\nu>n$,
this is a polynomial of degree~$n$ with non-vanishing
leading term
	$$\biggl(\Adach(N)\,{k^n\over2^nn!}\,f^*c^n\biggr)[N]
	={k^n\over2^nn!}\,\deg_{\Adach}f\cdot c^n[M]\ne 0\;.$$
From here, the proof proceeds as above.
\qed\enddemo

We will now concentrate on~$\C P^n$ and~$Q^n$
and prove \ZweitesGromovTheorem.
We have to show that a metric on~$\C P^n$ or on~$Q^n$
with non-negative Ricci curvature
is area-extremal in the sense of Gromov
(\Gromov, section~$5{4\over9}$):

\demo{Proof\/ \rm of \ZweitesGromovTheorem}
We start with~$M=\C P^n$.
Let~$H\to\C P^n$ be the hyperplane bundle,
then~$c:=c_1(H)$ generates~$H^2(M,\Z)$,
and the canonical line bundle of the complex structure on~$M$
is~$H^{-n-1}$.
We have~$w_2(N)=f^*w_2(M)=(n+1)\,f^*c\bmod 2$
for a spin map~$f\colon N\to M$,
and~$c^n[M]\ne0$.
There is a unique unitary connection on~$H$ such that the induced
connection on~$H^{-n-1}$ coincides with the connection on~$K$
which is induced by the Levi-Civita connection of the given K\"ahler metric~$g$.
Since~$\Ric_M\ge 0$, we have
	$$\kappa=\tr\(\Ric_M\)=2\norm{F^K}_g=2(n+1)\norm{\alpha}_g$$
with~$\alpha:=-iF^H$.
Thus, if~$f$ is area-nonincreasing and of non-vanishing $\Adach$-degree,
the estimate for~$\kapquer$ of \ZweitesGromovTheorem\ follows from
\NormLemma\ and \AdachDegreeCorollary\ with~$\alquer=f^*\alpha$.

Assume that~$\kapquer\ge\kappa\circ f>0$.
From the proof of \MainTheorem, we know that~$P_{f^*H}(k)$
vanishes precisely for~$k=n-1$, $n-3$, \dots, $1-n$.
This implies that the complex spinor bundle associated to the \spinc\ structure
with canonical bundle~$f^*H^{\pm(n+1)}$ admits a parallel spinor.
From here on,
the argument continues as in the proof of \GromovRigidity.

For~$M=Q^n$, we choose~$H$ to be the pull-back of the hyperplane bundle
on~$\C P^{n+1}$ via the canonical embedding.
Then~$P_{f^*H}(k)$ vanishes for~$k=n-2$, $n-4$, \dots, $2-n$.
Because~$P_{f^*H}(-k)=(-1)^n\,P_{f^*H}(k)$
and~$P_{f^*H}(k)$ has at most~$n$ zeros,
it follows that~$P_{f^*H}(n)=(-1)^nP_{f^*H}(-n)\ne 0$.
The rest of the proof is the same as above.
\qed\enddemo

\Remark\AdachProofRemark
One of the reasons that the Ricci curvature has to be strictly positive
in \GromovRigidity\ is to prevent the following pathological situation:
Suppose that~$i\colon M_0\to M$ is the inclusion of a complex submanifold
such that~$\norm{i^*F^K}_{i^*g}=\norm{F^K}_g$,
then the image of~$N$ could be contained in~$M_0$,
and~$N$ might not split as in the theorem.
In the setting of the theorem above,
this is prevented by the non-vanishing of the $\Adach$-degree of~$f$,
even if we allow~$M$ to have non-negative Ricci curvature
(the non-vanishing of the \spinc\ degree does not allow such a conclusion).
It would thus be interesting to know if the rigidity statement
of \ZweitesGromovTheorem\
continues to hold for K\"ahler metrics of non-negative Ricci curvature
on~$\C P^n$ or~$Q^n$.
\endremark

\Remark\GeneralHermitianSymmetricRemark
The proof above relies on the existence of a $k$-th root
of the canonical bundle of the K\"ahler metric on~$M$,
with~$k\ge n$.
Thus by the theorem of Kobayashi and Ochiai (\KO),
$M$ has to be either~$\C P^n$ or~$Q^n$,
and the argument above
has no obvious generalization to other compact Hermitian symmetric spaces.
Nevertheless, in view of \SpinCKapitel\ and the results of~\MinOo\ and~\GS,
one might hope that the statement of the theorem remains true for
all Ricci-positive K\"ahler metrics on compact Hermitian symmetric spaces.
\endremark

\Remark\ManyMetricsZwei
By the same argument as in \ManyMetrics,
there is an infinite-dimensional family of Ricci-positive K\"ahler metrics
on~$\C P^n$ and~$Q^n$ with the standard complex structure.
All these metrics are area-extremal in the sense of Gromov.
This suggests the following question:
does such an infinite-parameter family of area-extremal metrics
occur whenever~$(M,g)$ is area-extremal in the sense of Gromov
and the curvature of~$M$ is positive in a suitable sense
(e.g.\ Ricci-positive in the case of K\"ahler metrics)?
\endremark

\Section\IndexKapitel=4
\head 4. Estimates for Smooth Projective Varieties\endhead
\MainTheorem\ is particularly well adapted to symplectic manifolds
where the symplectic form~$\omega$ can be represented as the curvature
of a line bundle~$L$.
We now regard the following special case:
Suppose that~$M\subset\C P^N$ is a smooth algebraic variety,
or more generally, a symplectic submanifold,
equipped with an arbitrary Riemannian metric.
Let~$f\colon M\to\C P^n$ be an area-nonincreasing map that is homotopic
to the inclusion.
Then we get a rough estimate for the minimum
of the scalar curvature~$\kappa$ of~$M$.

If we know more about the topology of~$M$,
we can use the parameter~$k$ in the proof of \MainTheorem\
to ``fine-tune'' the estimate mentioned above.
For example, if~$M$ is a complete intersection,
we can determine the zeros of the Hilbert polynomial
to obtain a smaller upper bound for~$\min\kappa$.
This estimate will be sharp for all complete intersections of sufficiently
small total degree that carry a K\"ahler metric of constant positive
scalar curvature.

\Corollary\MainCorollary.
Let~$(M^{2n},g)$ be a compact, orientable Riemannian manifold,
and let~$f\colon M\to\C P^N$ be a smooth map,
such that~$\norm{f^*\omega}_g\le n$,
where~$\omega$ is the K\"ahler form of the Fubini-Study metric on~$\C P^N$
with constant holomorphic sectional curvature~$4$.
Assume that~$f^*[\omega]^n\ne 0\in H^{2n}(M,\R)$,
and that the second Stiefel-Whitney class~$w_2(M)$ of~$M$
is a multiple of~$f^*w_2(H)\in H^2(M,\Z_2)$,
where~$H$ is the hyperplane bundle of~$\C P^N$.
Then
	$$\kappa
	\not>\cases
		4n^2	&\text{if~$w_2(M)=n\,f^*w_2(H)$,\quad and}\\
		4n(n+1)	&\text{if~$w_2(M)=(n+1)\,f^*w_2(H)$.}
	\endcases$$
Assume that~$\kappa\ge4n^2$ in the first case
or~$\kappa\ge 4n(n+1)$ in the second case,
and that~$M$ is connected.
Then~$M$ is isometric to the complex quadric~$Q^n$ in the first case,
and to~$\C P^n$ in the second.

By \NormLemma,
the estimates above are applicable if~$f$ is an area-nonincreasing map:
In this case,
the absolute values of the~$n$ eigenvalues
of~$f^*\omega$ with respect to the induced metric
are all smaller than~$1$,
so~$\norm{f^*\omega}_g\le n$.

\demo{Proof\/ \rm of \MainCorollary}
We have equipped~$\C P^N$ with the Fubini-Study metric
of constant holomorphic sectional curvature~$4$,
i.e.\ scalar curvature~$4N(N+1)$.
In this case, the curvature of the hyperplane bundle
is given by~$F^H=-2i\omega$.
The image of the first Chern class~$c=c_1(H)$ of~$H$ in~$H^2(M,\R)$
is represented by~${1\over\pi}\omega$.
The estimate in the corollary follows from \MainTheorem\
with~$\alquer=2f^*\omega$,
because~$\norm{\smash{f^*\omega}}_g\le n$.

If we have~$\kappa\ge4n^2$ in the first case
or~$\kappa\ge 4n(n+1)$ in the second case,
then $M$ carries a metric of constant scalar curvature
and is biholomorphic to~$Q^n$ or~$\C P^n$ by \MainTheorem.
In particular, the group of complex automorphisms~$\frak A(M)$ of~$M$
acts transitively on~$M$.
From a theorem of Lichnerowicz (\Besse, Proposition~2.151),
it follows that a maximal compact connected subgroup of~$\frak A(M)$
acts by isometries.
It is easy to see that the metric on~$M$ is Hermitian symmetric,
and the rigidity statement of the corollary follows.
\qed\enddemo

\Rem
The topological conditions in \MainCorollary\ may be difficult to check
for arbitrary manifolds~$M$ and arbitrary maps~$f$.
There is however a large class of examples where these conditions
are automatically satisfied.
To begin with, regard~$\C P^N$ as a symplectic manifold.
If~$f^*\omega$ defines a symplectic structure on~$M$,
then~$f^*[\omega]^n$ is nonzero by definition
(actually it suffices that~$f$ be homotopic to a map~$\fquer$
for which~$\fquer^*\omega$ defines a symplectic structure).
For example, each smooth projective variety~$V\subset\C P^N$
is a symplectic submanifold.
Thus, $f^*[\omega]^n\ne 0$ whenever~$f\colon V\to\C P^N$
is homotopic to the inclusion.

Next, recall that~$M$ is spin iff~$w_2(M)$ vanishes,
and that~$f$ is a spin map iff~$w_2(M)=f^*w_2(\C P^N)$.
Because~$w_2(\C P^N)=(N+1)\,w_2(H)$,
we may conclude in both cases that~$w_2(M)$ is a multiple of~$f^*w_2(H)$.
\endremark

Now, we specialize the methods developped above
to study a particular type of smooth algebraic varieties:
A {\em complete intersection\/}~$V=V^n(a_1,\dots,a_r)$
of complex dimension~$n$
is the intersection of $r$ nonsingular hypersurfaces in~$\C P^{n+r}$
in general position,
defined by homogeneous polynomials of degrees~$a_1$, \dots, $a_r$
(cf.\ \Hirz, section~22.1).
We call~$\abs a:=a_1+\dots+a_r$ the {\em total degree\/} of~$V$.
We will look at arbitrary Riemannian metrics on~$V$
such that there exists a smooth map~$f\colon V\to\C P^{n+r}$
which is homotopic to the identity and area-nonincreasing
(actually, $\norm{f^*\omega}_g\le n$ is sufficient).
In this case,
B\"ar and Bleecker conjectured that~$\min\kappa\le 4n(n+1)$.
This bound has already been established in \MainCorollary.

In \IntroComplIntTheorem,
we have stated a stronger estimate
which fits well with a calculation of the average scalar curvature
of~$V$ with respect to the induced K\"ahler metric on~$V$ (\Ogiue).
Since we have explicit formulas for~$w_2(V)$ and~$\Adach(V)$,
we can find the minimal~$k$ such that the index of the operator~$D_k$
constructed in the proof of \MainTheorem\ above does not vanish.
We restate \IntroComplIntTheorem\ in a slightly stronger version:

\Theorem\ComplIntTheorem.
Let~$V=V^n(a_1,\dots,a_r)$ be a complete intersection,
equipped with an arbitrary Riemannian metric~$g$,
and set~$\abs a:=a_1+\dots+a_r$.
Let~$f\colon V\to\C P^{n+r}$ be homotopic to the inclusion,
such that~$\norm{f^*\omega}_g\le n$.
Then
	$$\min_{p\in V}\kappa(p)
	\le\cases
		4n\(n+r+1-\abs a\)
			&\text{if~$\abs a\le n+r$,}\\
		0	&\text{if~$\abs a>n+r$, $n$ is even,
				and~$V$ is spin,}\\
		4n	&\text{if~$\abs a>n+r$,
				and~$V$ is not spin,\quad and}\\
		8n	&\text{if~$\abs a>n+r$, $n$ is odd, and~$V$ is spin.}
	\endcases$$
\parr
If~$V$ is connected, $\abs a\le n+r$ (or~$\abs a=n+r+1$ and~$n$ is even)
and~$\kappa\ge 4n(n+r+1-\abs a)$,
then~$g$ is a K\"ahler metric
of constant scalar curvature~$4n(n+r+1-\abs a)$.
If moreover, $f$ is area-nonincreasing,
then~$V$ is a K\"ahler-Einstein manifold,
and~$f$ is an isometric, holomorphic immersion.

Before we prove \ComplIntTheorem,
let us look at some cases where our estimate is sharp.
Clearly, no upper bound for~$\min\kappa$ in a theorem as the above
can be negative:
let~$g$ be an arbitrary Riemannian metric.
Rescaling by some very large constant~$C\gg 0$,
we can make the absolute value of~$\min\kappa$ as small as we want
without violating the assumption that~$\norm{f^*\omega}\le n$.
In particular, $0$ is always the best estimate we can hope for.

\Remark\OgiueHanoRemark
By a computation of Ogiue,
the average scalar curvature of a complete intersection,
equipped with the pull-back of the Fubini-Study metric,
is precisely~$4n(n+r+1-\abs a)$ (\Ogiue, cf.~\BB).
For~$\abs a\le n+r$,
this value coincides with our estimate for~$\min\kappa$.
In this case,
\ComplIntTheorem\ gives a strong rigidity statement
if~$\kappa\ge 4n(n+r+1-\abs a)$:
It says that~$V$ carries a K\"ahler metric with constant scalar curvature
(which is not necessarily induced from the Fubini-Study metric).
According to \Hano,
the only complete intersections where the induced metric has constant
scalar curvature
are precisely the complex projective space~$\C P^n$ and the quadric~$Q^n$.
On the other hand,
some complete intersections admit K\"ahler-Einstein metrics which
are not induced from the Fubini-Study metric.
For example, Tian showed that the Fermat hypersurfaces
	$$V=\bigl\{\,[z_0:\ldots:z_{n+1}]\in\C P^{n+1}
		\bigm|z_0^a+\dots+z_{n+1}^a=0\,\bigr\}$$
admit such a metric if~$a=n$ or~$a=n+1$ (\Tian).
The authors do not know if there is an immersion~$f\colon V\to\C P^N$
for some~$N>n$
that is homotopic to the natural inclusion map
and satisfies~$\norm{f^*\omega}_g\le n$.
\endremark

Finally, as pointed out in~\BB,
any estimate for~$\min\kappa$ which is based on~$n$ and~$\abs a-r$ alone
must be positive if~$V$ is not spin:
The reason is that~$V=V^n(a)$ is a simply connected hypersurface if~$n\ge 3$,
which is not spin if~$n-a$ is odd.
In this case,
$V$ carries a metric of positive scalar curvature by a theorem
of Gromov and Lawson (\GL),
even if~$a$ is very large.
This explains the positive estimate in the third case of \ComplIntTheorem.
Next, choose~$n=4k+3$ and~$a$ odd,
then~$V^n(a)$ is spin and simply connected.
The generalized $\Adach$-genus~$\alpha(M)$ vanishes
in real dimension~$8k+6$,
so~$V$ carries a metric of positive scalar curvature
by a theorem of Stolz (\StEins),
regardless of the size of~$a$.
Thus, in the last case of \ComplIntTheorem,
the estimate again must be positive.
Nevertheless,
since our methods produce upper bounds for~$\min\kappa$
which are multiples of~$4n$,
one may expect that the estimates for the last two cases
are rather weak.

Before we prove \ComplIntTheorem,
we recall a few facts about complete intersections (\Hirz):

\Lemma\ComplIntLemma.
Let~$V=V^n(a_1,\dots,a_r)$ be a complete intersection,
and let~$f\colon V\to\C P^{n+r}$ be homotopic to the inclusion.
Let~$\bold x:=f^*c_1(H)\in H^2(V,\Z)$.
Then
$$\align
	w_2(V)&\equiv\(n+r+1-\abs a\)\,\bold x\quad\pmod 2\;,\\
		\text{and}\qquad
	\Adach(V)
	&=\biggl({\bold x/2\over\sinh(\bold x/2)}\biggr)^{n+r+1}
		\,\prod_{j=1}^r{\sinh(a_j\bold x/2)\over a_j\bold x/2}\;.
\endalign$$
In particular,
$V$ is spin iff~$\abs a\equiv n+r+1\pmod2$.
Moreover,
	$$\bold x^n[V]=a_1\cdots a_r\;.$$

\Proof
Clearly, the cohomology class~$\bold x$ depends only on the homotopy class
of the map~$f$.
Now, the first two equations follow from the fact
that the total Chern class of~$TV$ is given by
	$$c(V)=(1+\bold x)^{n+r+1}\,\prod_{j=1}^r\(1+a_j\bold x\)^{-1}\;,$$
cf.~\Hirz, chapter~22, equation~(1).
The last equation is also contained in~\Hirz.
\qed\enddemo

\demo{Proof\/ \rm of \ComplIntTheorem}
On~$V$, there exists a \spinc\ structure~$P_\Spinc^k$
with canonical line bundle~$f^*H^k$
iff~$k\,f^*w_2(H)=w_2(V)\in H^2(V,\Z_2)$.
By~\ComplIntLemma,
this is the case if~$k\equiv n+r+1-\abs a\pmod2$.
We denote the Dirac operator associated to~$P_\Spinc^k$ by~$D_k^V$.
We will determine the smallest~$k\ge 0$ such that~$D^V_k$ exists
and~$\ind(D^V_k)\ne 0$.
Then our assertion will follow from \Estimate\
as in the proof of \MainTheorem.

By \ComplIntLemma,
the index of~$D^V_k$ can be calculated as follows:
\NeueFormel\IndexRechnung
$$\multline
	\ind\(D_k^V\)
	=\Bigl(\Adach(V)\,e^{\textstyle{k\bold x\over2}}\Bigr)[V]
	=\Biggl(\biggl({\bold x/2\over\sinh(\bold x/2)}\biggr)^{n+r+1}
		\,e^{\textstyle{k\bold x\over2}}
		\,\prod_{j=1}^r{\sinh(a_j\bold x/2)\over a_j\bold x/2}
		\Biggr)[V]
\allowdisplaybreak
\split
	&=2^{-n-1}\,\res_{x=0}\biggl(\sinh\Bigl({x\over2}\Bigr)^{-n-r-1}
		\,e^{\textstyle{kx\over2}}
		\,\prod_{j=1}^r\sinh{a_jx\over2}\biggr)
\allowdisplaybreak
	&=2^{-n-1}\,\res_{x=0}\Biggl(\sinh\Bigl({x\over2}\Bigr)^{-n-1}
		\,e^{\textstyle{kx\over2}}
		\,\prod_{j=1}^r{e^{a_jx\over2}-e^{-{a_jx\over2}}
				\over e^{x\over2}-e^{-{x\over2}}}\Biggr)
\allowdisplaybreak
	&=2^{-n-1}\,\res_{x=0}\Biggl(\sinh\Bigl({x\over2}\Bigr)^{-n-1}
		\,\sum_{l_1={1-a_1\over2}}^{a_1-1\over2}\,\cdots
		\,\sum_{l_r={1-a_r\over2}}^{a_r-1\over2}
		\,e^{\textstyle{(2l_1+\dots+2l_r+k)\,x\over2}}\Biggr)
\allowdisplaybreak
	&=\sum_{l_1={1-a_1\over2}}^{{a_1-1\over2}}\,\cdots
		\,\sum_{l_r={1-a_r\over2}}^{{a_r-1\over2}}
		P_n\(2l_1+\dots+2l_r+k\)\;.
\endsplit\endmultline\tag\FormelNummer\IndexRechnung$$
Here, $P_n(k')$ denotes the index of the operator~$D_{k'}^{\C P^n}$
constructed on~$\C P^n=V^n(1,\dots,1)\subset\C P^{n+r}$.
Because in the multiple sum above,
	$$2l_1+\dots+2l_r+k\equiv\abs a-r+k\equiv n+1\pmod 2\;,$$
the operators~$D_{2l_1+\dots+2l_r+k}^{\C P^n}$ exist indeed.

Let us regard the case that~$n$ is even.
In this case, $P_n(k')$ is an even function of~$k'$,
and~$P_n(k')>0$ whenever~$\abs{k'}>n-1$ and~$k'\equiv n+1\pmod 2$
by~\CPnIndexEstimate.
The multiple sum in~\IndexRechnung\ contains terms of the form~$P_n(k')$,
where~$k'$ ranges between~$r-\abs a+k$ and~$\abs a-r+k$.
If~$\abs a\le n+r$,
then the smallest~$k\in\Z$ with~$k\equiv n+r+1-\abs a$
such that~$P_n(\abs a-r+k)>0$ is given by
	$$k_0=n+r+1-\abs a\;.$$
By \Estimate, we have in this case
	$$\min_{p\in V}\kappa(p)\le 4n\(n+r+1-\abs a\)\;.$$
On the other hand, if~$\abs a>n+r$,
we may choose~$k=0$ if~$V$ is spin and~$k=1$ otherwise.
The corresponding scalar curvature estimates are
	$$\min_{p\in V}\kappa(p)\le 0\quad\text{if~$V$ is spin,\qquad and}
	\qquad\min_{p\in V}\kappa(p)\le 4n\quad\text{otherwise.}$$
This proves the estimates of \ComplIntTheorem\ if~$n$ is even.

If~$n$ is odd, then~$P_n(k')$ is an odd function in~$k'$.
In particular by~\IndexRechnung,
$\ind(D_0^V)=0$ if~$V$ is spin,
even if~$\abs a$ is very large.
We calculate
$$\multline
	\ind\(D_k^V\)-\ind\(D_{k-2}^V\)
	=2^{-n-1}\,\res_{x=0}\biggl(\sinh\Bigl({x\over2}\Bigr)^{-n-r-1}
		\,\Bigl(e^{\textstyle{kx\over2}}
			-e^{\textstyle{(k-2)x\over2}}\Bigr)
		\,\prod_{j=1}^r\sinh{a_jx\over2}\biggr)\\
	=2^{-n}\,\res_{x=0}\biggl(\sinh\Bigl({x\over2}\Bigr)^{-n-r}
		\,e^{\textstyle{(k-1)x\over2}}
		\,\prod_{j=1}^r\sinh{a_jx\over2}\biggr)
	=\ind\(D_{k-1}^W\)\;,
\endmultline$$
where~$W:=V^{n-1}(1,a_1,\dots,a_r)$ is the transverse intersection of~$V$
with a generic hyperplane in~$\C P^{n+r}$
(actually, we do not need the existence of~$W$.
We only need the properties of the formal expression for~$\ind(D_{k-1}^W)$,
which we established in the previous paragraphs).
Let~$k_0\ge0$ be the smallest value of~$k$
such that~$\ind(D_k^W)\ne0$.
It follows that~$k_0+1$ is the smallest value of~$k$
such that~$\ind(D_k^V)\ne0$.
This completes the proof of the estimates.

Let us now look at the case where~$\kappa\ge 4n(n+r+1-\abs a)$
and~$\abs a\le n+r$ (or~$\abs a=n+r+1$ and~$n$ is even),
Then we have
	$$\kappa\ge 4n(n+r+1-\abs a)\ge
	4\,(n+r+1-\abs a)\,\norm{f^*\omega}_g
	=2\norm{f^*F^{H^{-(n+r+1-\abs a)}}}_g\;,$$
and the corresponding Dirac operator has non-vanishing index
by our calculations above.
Then as in the proof of \Rigidity,
there exists a parallel spinor~$\psi$ on~$V$.
Since we assumed that~$f$ is homotopic to the inclusion~$V\subset\C P^N$,
we know that~$\(f^*\omega^n\)[V]\ne 0$,
in particular~$f^*\omega^n\ne 0$ in at least one point of~$V$.
This is enough to ensure that~$\psi$ is a parallel, pure spinor,
thus it defines a parallel complex structure on~$V$,
and~$g$ is a K\"ahler metric of constant scalar curvature.

If we assume in addition that~$f$ is area-nonincreasing,
then each of the $n$ eigenvalues of~$f^*\omega$ has to be~$1$,
so~$f$ is an isometric immersion,
and~$f_*T_pV$ is a complex subspace of~$T_{f(p)}\C P^{n+r}$
for all~$p\in V$.
Moreover, the canonical line bundle of~$V$ has curvature
	$$-2i\,(n+r+1-\abs a)\,f^*\omega=-(n+r+1-\abs a)\,f^*F^H\;.$$
Because~$V$ has positive Ricci curvature by \Rigidity,
this implies that~$f$ is holomorphic
if we choose the proper complex structure on~$V$.
Finally, we conclude that the metric~$g$ is K\"ahler-Einstein
with Einstein constant~$2(n+r+1-\abs a)$.
\qed\enddemo

\Refs
\widestnumber\key{KW}

\Quelle\BB[BB]\Preprint
  C.~B\"ar, D.~Bleecker:
  The Dirac Operator and the Scalar Curvature of Continuously Deformed
  Algebraic Varieties,
  preprint\miscnote to appear in Contemporary Mathematics

\Quelle\Besse[Be]\Buch
  A.~L.~Besse:
  Einstein Manifolds,
  Spinger,  Berlin-Heidel\-berg-New York (1987)

\Quelle\GS[GS]\Preprint
  S.~Goette, U.~Semmelmann:
  Scalar Curvature Estimates for Compact Symmetric Spaces,
  \miscnote in preparation

\Quelle\Gromov[G]\ImBuch
  M.~Gromov:
  {Positive curvature, macroscopic dimension, spectral gaps
  and higher signatures},
  in S.~Gindikin, J.~Lepowski, R.~L.~Wilson:
  {Functional Analysis on the Eve of the 21st Century, Vol.~II},
  Progress in Mathematics Vol.~132 (1996), 1--213

\Quelle\GL[GL]\Journal
  M.~Gromov, H.~B.~Lawson Jr.:
  The classification of simply-connected manifolds
  of positive scalar curvature,
  Ann.\ of Math.~111 (1980), 423--434

\Quelle\Hano[Ha]\Journal
  J.~Hano:
  Einstein complete intersections in complex projective spaces,
  Math.\ Ann.~216 (1975), 197--208

\Quelle\Hirz[Hir]\Buch
  F.~Hirzebruch:
  Topological Methods in Algebraic Geometry,
  Third Enlarged Edition, Springer, New York (1966)

\Quelle\Hitchin[Hit]\Journal
  N.~Hitchin:
  Harmonic Spinors,
  Adv.\ in Math.~14 (1974), 1--55

\Quelle\KW[KW]\Journal
  J.~L.~Kazdan, F.~W.~Warner:
  Scalar curvature and conformal deformation of Riemannian structure,
  J.\ Diff.\ Geom.~10 (1975), 113--134

\Quelle\KO[KO]\Journal
  S.~Kobayashi, T.~Ochiai:
  Characterizations of complex projective spaces and hyperquadrics,
  J.\ Math.\ Kyoto Univ.~13 (1973), 31--47

\Quelle\Kramer[Kr]\Preprint
  W.~Kramer:
  Der Dirac-Operator auf Faserungen,
  Dissertation, Bonner Mathematische Schriften~317,
  Universit\"at Bonn~1999

\Quelle\LM[LM]\Buch
  H.~B.~Lawson Jr., M.-L.~Michelsohn:
  Spin Geometry,
  Princeton Univ.\ Press, Princeton, N.\ J.\ (1989)

\Quelle\LlaZwei[L1]\Journal
  M.~Llarull:
  Scalar curvature estimates for $(n+4k)$-dimensional manifolds,
  Diff.\ Geom.\ Appl.~6 (1996), 321--326

\Quelle\LlaEins[L2]\Journal
  M.~Llarull:
  Sharp Estimates and the Dirac Operator,
  Math.\ Ann.~310 (1998), 55--71

\Quelle\MinOo[M]\ImBuch
M.~Min-Oo:
Scalar Curvature Rigidity of Certain Symmetric Spaces,
in :
{Geometry, Topology and Dynamics}\procinfo{Montreal, PQ, 1995},
CRM Proc.\ Lecture Notes, 15,
Amer.\ Math.\ Soc., Providence, RI (1998), 127--136

\Quelle\Moro[Mo]\Journal
  A.~Moroianu:
  Parallel and Killing Spinors on \SPinc\ Manifolds,
  Commun.\ Math.\ Phys.~187 (1997), 417--427

\Quelle\Ogiue[O]\Journal
  K.~Ogiue:
  Scalar curvature of complex submanifolds in complex projective spaces,
  J.\ Diff.\ Geom.~5 (1971), 229--232

\Quelle\Schro[Sch]\Journal
  E.~Schr\"odinger:
  Diracsches Elektron im Schwerefeld.~I,
  Sitzungsber.\ Preuss.\ Akad.\ Wiss.~11 (1932), 105--128

\Quelle\StEins[S1]\Journal
  S.~Stolz:
  Simply connected manifolds of positive scalar curvature,
  Bull.\ Am.\ Math.\ Soc.~23 (1990), 427-432

\Quelle\StZwei[S2]\ImBuch
  S.~Stolz:
  Positive scalar curvature metrics---existence and classification questions,
  in : {Proc.\ of the ICM, Vol.~I}%
  \procinfo{ICM Z\"urich, 1994},
  Birkh\"auser, Basel (1995), 625--636

\Quelle\Tian[T]\Journal
  G.~Tian:
  On K\"ahler-Einstein Metrics on Certain Manifolds with~$c_1(M)>0$,
  Inv.\ Math.~89 (1987), 225--246

\endRefs

\enddocument